\documentclass[12pt]{amsart}
\usepackage[francais]{babel}
\usepackage{vmargin}
\usepackage[dvips]{graphicx}
\usepackage{color}
\input{amssym.def}
\input{amssym}
\usepackage{a4wide}
\usepackage{supertabular}
\usepackage{array}
\usepackage{multicol}

\setmarginsrb{3cm}{2cm}{3cm}{3cm}{75pt}{20pt}{20pt}{30mm}
\def\derpar#1#2{\frac{\partial#1}{\partial#2}}
\def\R{\mathbb R}
\def\N{\mathbb N}
\def\T{\mathbb T}
\def\C{\mathbb C}

\def\Z{\mathbb Z}

\def\var{\varepsilon}

\def\pa{\partial}

\def\Om{\Omega}

\def\ov{\overline}
\def\cal{\mathcal}

\def\hat{\widehat}

\def\tilde{\widetilde}

\def\Id{{\rm Id}\,}
\def\dps{\displaystyle}

\def\<{\langle}
\def\>{\rangle}

\def\cF{{\cal F}} 

\def\cZ{{\cal Z}}

\def\med{\medskip}
\def\sm{\smallskip}
\def\bul{$\bullet$\ }

\def\begeq{\begin{equation}}
\def\endeq{\end{equation}}
\def\begar{\begin{eqnarray}}
\def\endar{\end{eqnarray}}
\def\begar*{\begin{eqnarray*}}
\def\endar*{\end{eqnarray*}}
\def\begal{\begin{align}}
\def\endal{\end{align}}
\def\begal*{\begin{align*}}
\def\endal*{\end{align*}}

\newtheorem{Thm}{Th\'eor\`eme}

\theoremstyle{definition}

\theoremstyle{remark}

\newtheorem*{Thm*}{Theorem}
\newtheorem*{Lem*}{Lemma}
\newtheorem*{Conj*}{Conjecture}
\newtheorem*{Cor*}{Corollary}
\newtheorem*{Def*}{Definition}
\newtheorem*{Prop*}{Proposition}
\newtheorem*{Exo*}{Exercise}
\newtheorem*{Exs*}{Examples}
\newtheorem*{Ex*}{Example}
\newtheorem*{Rk*}{Remark}
\newtheorem*{Rks*}{Remarks}

\def\signcv{\bigskip \begin{center} {\sc C\'edric Villani\par\vspace{2mm}
ENS Lyon \& Institut Universitaire de France\par
UMPA, UMR CNRS 5669\par
46 all\'ee d'Italie\par
69364 Lyon Cedex 07\par
FRANCE\par\vspace{2mm}
e-mail:} \tt{cvillani@umpa.ens-lyon.fr} \end{center}}

\def\signcm{\bigskip \begin{center} {\sc Cl\'ement
      Mouhot\par\vspace{2mm}
DAMTP, University of Cambridge,\par
{\it On leave from:} 
CNRS \& \'Ecole Normale Sup\'erieure\par
DMA, UMR CNRS 8553\par
45 rue d'Ulm, F 75320 Paris Cedex 05\par
FRANCE\par\vspace{2mm}
e-mail:} \tt{Clement.Mouhot@ens.fr} \end{center}}

\begin{document}

\title{Amortissement Landau}

\author{C. Mouhot}
\author{C. Villani}

\maketitle

\vspace*{-10mm}

\begin{abstract} 
  Dans cette note nous pr\'esentons les principaux r\'esultats du
  r\'ecent travail \cite{MV}, o\`u le ph\'enom\`ene d'amortissement
  Landau est pour la premi\`ere fois \'etabli dans un contexte non
  lin\'eaire.  \med

  \noindent {\sc Abstract.} In this note we present the main results
  from the recent work \cite{MV}, which for the first time establish
  Landau damping in a nonlinear context.
\end{abstract}

\med

{\bf Mots-cl\'es.} amortissement Landau; physique des plasmas;
astrophysique; \'equation de Vlasov--Poisson.

\section{Introduction}
\label{sec:intro}

Le \og mod\`ele standard\fg{} de la physique des plasmas classique est
l'\'equation de Vlasov--Poisson--Landau \cite{landau,vlasov}, ici
\'ecrite avec des conditions aux limites p\'eriodiques et en unit\'es
adimensionn\'ees: 
\begeq\label{VL} \derpar{f}{t} + v\cdot\nabla_x f +
F[f]\cdot\nabla_v f = \frac{\log\Lambda}{2\pi \Lambda}\, Q_L(f,f),
\endeq 
o\`u $f=f(t,x,v)$ est la fonction de distribution des \'electrons
($t\geq 0$, $v\in\R^3$, $x\in\T^3 = \R^3/\Z^3$), \begeq\label{Ff}
F[f](t,x) = -\iint \nabla W(x-y)\, f(t,y,w)\,dw\,dy
\endeq
est la force auto-induite, $W(x) = 1/|x|$ est le potentiel
d'interaction coulombien, et $Q_L$ est l'op\'erateur de collision de
Landau, d\'ecrit par exemple dans \cite{LL:kin:81} ou \cite{vill:handbook:02}.  Le
param\`etre $\Lambda$ est tr\`es grand, variant typiquement entre
$10^2$ et $10^{30}$. 

Sur de tr\`es grandes \'echelles de temps (disons $O(\Lambda/\log
\Lambda)$), les ph\'enom\`enes dissipatifs jouent un r\^ole non
n\'egligeable, et l'augmentation de l'entropie est suppos\'ee forcer
la convergence (lente) vers une maxwellienne. Gr\^ace aux progr\`es
r\'ecents sur l'hypocoercivit\'e, ce m\'ecanisme est maintenant assez
bien compris math\'ematiquement parlant, d\`es que l'on dispose d'estimations de r\'egularit\'e
globale (voir \cite{vill:hypoco} et les r\'ef\'erences incluses).

Dix ans apr\`es avoir mis au point ce sc\'enario collisionnel, Landau
\cite{landau} formulait une pr\'ediction beaucoup plus subtile: la
stabilit\'e d'un \'equilibre homog\`ene v\'erifiant certaines
conditions --- par exemple fonction radiale de la vitesse $v$, mais
pas forc\'ement gaussienne --- sur des \'echelles de temps beaucoup
plus courtes (disons $O(1)$), par le jeu de m\'ecanismes purement
conservatifs.  Ce ph\'enom\`ene, appel\'e {\bf amortissement Landau},
est une propri\'et\'e de l'\'equation (non collisionnelle!) de Vlasov, obtenue
en posant $\Lambda=\infty$ dans \eqref{VL}. La d\'ecouverte de Landau
fut jug\'ee ``stup\'efiante'' par ses contemporains; c'est maintenant l'une des bases
th\'eorique de la physique classique des plasmas (parmi un tr\`es
grand nombre de r\'ef\'erences citons \cite{akhiezer,LL:kin:81}). Des effets
d'amortissement similaires ont \'et\'e pr\'edits en astrophysique \cite{LB:damping},
mais aussi dans d'autres domaines de la physique.

L'amortissement Landau est compris depuis longtemps au niveau de
l'\'equation lin\'earis\'ee \cite{backus,degond:landau,hayes:fluids,hayes:cimento,MF:landaudamping,saenz}. 
Cependant, il y a d\'ej\`a cinquante ans de cela, dans le premier volume du
{\em Journal of Mathematical Physics}, Backus \cite{backus} notait que l'\'echelle de
temps de la lin\'earisation est parfois bien plus courte que l'\'echelle de temps de l'amortissement
Landau, et exprimait des doutes sur la pertinence de l'approche par lin\'earisation.
Son objection n'a jamais \'et\'e r\'efut\'ee, car les difficult\'ees conceptuelles et techniques
associ\'ees \`a l'\'equation non lin\'eaire ont s\'ev\`erement limit\'e les r\'esultats obtenus dans
ce cadre: dans \cite{caglimaff:VP:98,HV:landau} on prouve seulement l'existence
de {\em certaines} solutions amorties. 

Nous comblons cette lacune dans un travail r\'ecent \cite{MV}, dont nous d\'ecrivons maintenant le
r\'esultat principal.

\section{R\'esultat principal}
\label{sec:resultat}

Si $f$ est une fonction d\'efinie sur $\T^d\times\R^d$, on note, pour
$k\in\Z^d$, $\eta\in\R^d$,
\[ \hat{f}(k,v) = \int_{\T^d} f(x,v)\, e^{-2i\pi k\cdot x}\,dx,\qquad
\tilde{f}(k,\eta) = \iint_{\T^d \times \R^d} 
   f(x,v)\,e^{-2i\pi k\cdot x}\,e^{-2i\pi \eta\cdot v}\,dv\,dx.\]
On pose en outre
\begin{equation}\label{norm1}
 \|f\|_{\lambda,\mu,\beta} = 
\sup_{k,\eta} 
\Bigl( |\tilde{f}(k,\eta)|\,e^{2\pi\lambda|\eta|}\, e^{2\pi \mu |k|}\Bigr)
+ \iint_{\T^d \times \R^d} |f(x,v)|\,e^{2\pi \beta|v|}\,dv\,dx.
\end{equation}

\begin{Thm}[amortissement Landau non lin\'eaire pour interactions
  g\'en\'erales]
\label{thmlim}
Soient $d\geq 1$, et $f^0:\R^d\to\R_+$ un profil de vitesses
analytique. Soit $W:\T^d\to\R$ un potentiel d'interaction. Pour
$k\in\Z^d$, $\xi\in\C$, on pose
\[ {\cal L}(k,\xi) = -4\pi^2\, \hat{W}(k) \int_0^\infty e^{2\pi
  |k|\xi^*t}\,|\tilde{f}^0(kt)|\,|k|^2\,t\,dt.\] On suppose qu'il
existe $\lambda>0$ tel que
\begeq\label{condf0} 
\sup_{\eta\in\R^d}\ |\tilde{f}^0(\eta)|\,
e^{2\pi\lambda|\eta|} \leq C_0,\qquad \sum_{n\in\N^d}
\frac{\lambda^n}{n!} \|\nabla_v^n f^0\|_{L^1(dv)} \leq C_0,
\endeq
\begeq\label{condL}
\inf_{k\in\Z^d}\ \inf_{0\leq \Re \, \xi \, < \lambda}\ 
\bigl|{\cal L}(k,\xi) - 1 \bigr| \geq\kappa>0,
\endeq
\begeq\label{condW}
\exists \, \gamma \ge 1;\ \forall \, k\in \Z^d;\qquad
|\hat{W}(k)| \leq \frac{C_W}{|k|^{1+\gamma}}.
\endeq
Alors d\`es que $0<\lambda'<\lambda$, $0<\mu'<\mu$, $\beta>0$,
$r\in\N$, il existe $\var>0$ et $C>0$, d\'ependant seulement de
$d,\gamma,\lambda,\lambda',\mu,\mu',C_0,\kappa,C_W,\beta,r$ tels que,
si $f_i\geq 0$ v\'erifie \begeq\label{condfi} \delta:=
\|f_i-f^0\|_{\lambda,\mu,\beta} \leq \var,
\endeq
alors l'unique solution de l'\'equation de Vlasov non lin\'eaire
\begeq\label{nlv}
\derpar{f}{t} + v\cdot\nabla_x f + F[f]\cdot\nabla_v f = 0,\qquad
F[f](t,x) = -\int \nabla W(x-y) \, f(t,y,w)\,dw\,dy,
\endeq
d\'efinie pour tous les temps et telle que $f(0,\,\cdot\,)=f_i$, v\'erifie
\begeq\label{rr}
\bigl\|\rho(t,\,\cdot\,)-\rho_\infty\bigr\|_{C^r(\T^d)} \leq C\,\delta\,e^{-2\pi \lambda'|t|},
\endeq
o\`u $\rho(t,x) = \int f(t,x,v)\,dv$, $\rho_\infty = \iint
f_i(x,v)\,dv\,dx$. En outre il existe des profils analytiques
$f_{+\infty}(v)$, $f_{-\infty}(v)$ tels que
\[ f(t,\,\cdot\,) \xrightarrow[]{t\to\pm\infty} f_{\pm\infty} \qquad
\text{faiblement}\]
\[ \int f(t,x,\,\cdot\,)\,dx \xrightarrow[]{t\to\pm\infty}
f_{\pm\infty}\qquad \text{fortement (dans $C^r(\R^d_v)$)},\] 
ces convergences \'etant \'egalement en $O(\delta\,e^{-2\pi\lambda'|t|})$.
\end{Thm}
\med

Ce th\'eor\`eme, enti\`erement constructif, est presque optimal comme
le montrent les commentaires qui suivent.  
\med

{\bf Commentaires sur les hypoth\`eses:} Les conditions aux limites
p\'eriodiques sont bien s\^ur contestables; cependant, au vu des
contre-exemples de Glassey et Schaeffer \cite{GS}, un m\'ecanisme de
confinement --- ou tout au moins de limitation des longueurs d'onde
spatiales --- est indispensable. La condition \eqref{condf0} exprime
quantitativement l'analyticit\'e du profil $f^0$, sans laquelle on ne
peut esp\'erer avoir de convergence exponentielle.  L'in\'egalit\'e
\eqref{condL} est une condition de stabilit\'e lin\'eaire,
essentiellement optimale, qui couvre les cas physiquement
int\'eressants~: l'interaction de Newton aux longueurs d'onde
inf\'erieures \`a la longueur d'instabilit\'e de Jeans, et
l'interaction coulombienne pour toutes les longueurs d'ondes autour de
profils $f^0$ radialement sym\'etriques en dimension sup\'erieure ou
\'egale \`a 3. La condition \eqref{condW} en revanche n'intervient que
dans la stabilit\'e non lin\'eaire, elle inclue les interactions de
Coulomb ou Newton dans le cas limite $\gamma =1$. Quant \`a la
condition \eqref{condfi}, son caract\`ere perturbatif est naturel au
vu des sp\'eculations th\'eoriques et des \'etudes num\'eriques sur le
sujet.  \med

{\bf Commentaires sur les conclusions}

\begin{enumerate}
\item La convergence en temps grand est bas\'ee sur un m\'ecanisme
  r\'eversible, purement d\'eterministe, sans fonctionnelle de
  Lyapunov ni interpr\'etation variationnelle. Les profils
  asymptotiques $f_{\pm\infty}$ gardent d'ailleurs la m\'emoire de la
  donn\'ee initiale et de l'interaction. Cette convergence \og sans
  raison d'\^etre\fg{} n'\'etait pas vraiment attendue, puisque la
  th\'eorie quasilin\'eaire de l'amortissement Landau \cite[Vol. II,
  Section 9.1.2]{akhiezer} ne pr\'edit la convergence qu'apr\`es une
  prise de moyenne sur des ensembles statistiques.

\item On peut interpr\'eter ce r\'esultat dans l'esprit du
  th\'eor\`eme KAM: pour l'\'equation de Vlasov lin\'eaire, la
  convergence est forc\'ee par une infinit\'e de lois de conservation
  qui rendent le mod\`ele \og compl\`etement int\'egrable\fg{}~; d\`es que
  l'on ajoute un couplage non lin\'eaire, les lois de conservation
  disparaissent mais la convergence demeure.



\item \'Etant donn\'e un \'equilibre stable $f^0$, un voisinage entier
  --- en topologie analytique --- de $f^0$ est rempli par des
  trajectoires homoclines ou (en g\'en\'eral) h\'et\'eroclines. C'est
  la dimension infinie qui permet ce comportement remarquable de
  l'\'equation de Vlasov non lin\'eaire.

\item La convergence en temps grand de la fonction de distribution n'a
  lieu qu'au sens faible; les normes des d\'eriv\'ees en vitesse
  croissent tr\`es vite en temps grand, ce qui traduit une
  filamentation dans l'espace des phases, et un transfert d'\'energie
  (ou d'information) des basses vers les hautes fr\'equences (\og
  turbulence  faible\fg{}).

\item C'est ce m\'ecanisme de transfert d'information aux petites
  \'echelles qui permet de r\'econcilier la r\'eversibilit\'e de
  l'\'equation de Vlasov--Poisson avec l'apparente irr\'eversibilit\'e
  de l'amortissement Landau. Notons que le m\'ecanisme \og dual\fg{} de
  transfert d'\'energie vers les grandes \'echelles, aussi appel\'e
  {\em radiation}, a \'et\'e largement \'etudi\'e dans le cadre des
  syst\`emes hamiltoniens de dimension infinie.

\end{enumerate}

On trouvera davantage de commentaires, aussi bien math\'ematiques
que physiques, dans \cite{MV}.

\section{Stabilit\'e lin\'eaire}
\label{sec:lin}

La stabilit\'e lin\'eaire est la premi\`ere \'etape de notre \'etude;
elle ne demande qu'un investissement technique assez r\'eduit.

Lin\'earis\'ee autour d'un \'equilibre homog\`ene $f^0(v)$,
l'\'equation de Vlasov devient 
\begeq\label{LV} \derpar{h}{t} +
v\cdot\nabla_x h - (\nabla W\ast \rho)\cdot\nabla_v f^0 =0,\qquad \rho
= \int h\,dv.
\endeq
Il est bien connu que cette \'equation se d\'ecouple en une infinit\'e
d'\'equations ind\'ependantes r\'egissant l'\'evolution des modes de
$\rho$: pour tout $k\in\Z^d$ et tout $t\geq 0$, 
\begeq\label{rho}
\hat{\rho}(t,k) - \int_0^t K^0(t-\tau,k)\,\hat{\rho}(\tau,k)\,d\tau =
\tilde{h}_i(k,kt),
\endeq
o\`u $h_i$ est la donn\'ee initiale, et $K^0$ un noyau int\'egral qui
d\'epend de $f^0$: 
\begeq\label{K0} K^0(t,k) = -
4\pi^2\,\hat{W}(k)\,\tilde{f}^0(kt)\,|k|^2\,t.
\endeq

On d\'eduit alors de r\'esultats classiques sur les \'equations de
Volterra que pour tout $k\neq 0$ la d\'ecroissance de
$\hat{\rho}(t,k)$ quand $t\to +\infty$ est essentiellement
contr\^ol\'ee par le pire de deux taux de convergence: \sm

\bul le taux de convergence du terme source au membre de droite de
\eqref{rho}, qui ne d\'epend que de la r\'egularit\'e de la donn\'ee
initiale dans la variable de vitesse; \sm

\bul $e^{-\lambda t}$, o\`u $\lambda$ est le plus grand r\'eel positif
tel que la transform\'ee de Fourier--Laplace (dans la variable $t$) de
$K^0$ ne s'approche pas de~1 dans la bande complexe $\{0\leq \Re z
\leq \lambda\}$.  \sm

Le probl\`eme consiste donc \`a trouver des conditions suffisantes sur
$f^0$ pour garantir la stricte positivit\'e de $\lambda$.  Depuis
Landau, cette \'etude est traditionnellement r\'ealis\'ee gr\^ace \`a
la formule d'inversion de la transform\'ee de Laplace; cependant, dans
la perspective de l'\'etude non lin\'eaire, nous lui pr\'ef\'erons une
approche plus \'el\'ementaire et constructive, bas\'ee sur la simple
formule d'inversion de Fourier.

On montre ainsi l'amortissement Landau lin\'eaire, sous les conditions
\eqref{condL} et \eqref{condf0}, pour n'importe quelle interaction $W$
telle que $\nabla W \in L^1(\T^d)$, et pour toute condition
initiale analytique (sans restriction de taille dans ce contexte
lin\'eaire).  On retrouve comme cas particulier tous les r\'esultats
pr\'ec\'edemment connus sur l'amortissement Landau lin\'eaire
\cite{degond:landau,MF:landaudamping,saenz}; mais on couvre
\'egalement l'interaction newtonienne. En effet, la condition
\eqref{condL} est v\'erifi\'ee d\`es que l'une ou l'autre des
conditions suivantes est vraie: \sm

(a) $\forall \, k\in\Z^d$, $\forall \, z\in\R$, $\hat{W}(k)\geq 0$, \
$z\,\phi'_k(z)\leq 0$, o\`u $\phi_k$ est la \og marginale\fg{} de $f^0$ selon
la direction $k$, d\'efinie par
\[ \phi_k(z) = \int_{\frac{kz}{|k|}+k^\bot} f^0(w)\,dw, \quad z \in \R
\ ;\]

(b) $\dps 4\pi^2\, \bigl(\max\,|\hat{W}(k)|\bigr)
\left(\sup_{|\sigma|=1} \int_0^\infty |\tilde{f}^0(r\sigma)|\,r\,dr
\right) < 1$.  \med

Cette derni\`ere condition couvre l'interaction de Newton en-de\c c\`a de la longueur
de Jeans. On renvoie \`a \cite[Section 3]{MV} pour plus de d\'etails.

\section{Stabilit\'e non lin\'eaire}
\label{sec:nonlin}

Pour traiter la stabilit\'e non lin\'eaire, nous commen\c cons par
introduire des {\bf normes analytiques} \og hybrides\fg{} (bas\'ees sur la taille
des d\'eriv\'ees successives dans la variable de vitesse, et sur la
taille des coefficients de Fourier dans la variable de position) et
\og glissantes\fg{} (la norme utilis\'ee changera au cours du temps pour
tenir compte des transferts dans l'espace des phases). Cinq indices
permettent d'obtenir toute la souplesse n\'ecessaire: \begeq\label{Z}
\|f\|_{\cZ^{\lambda,(\mu,\gamma);p}_\tau} = \sum_{k\in\Z^d}
\sum_{n\in\N^d} e^{2\pi\mu|k|}\, (1+|k|)^\gamma\,
\frac{\lambda^n}{n!}\, \Bigl\|\bigl(\nabla_v + 2i\pi \tau k\bigr)^n
\hat{f}(k,v)\Bigr\|_{L^p(dv)}.
\endeq
(Par d\'efaut $\gamma=0$.) Un th\'eor\`eme fastidieux d'injection
\og \`a la Sobolev\fg{} compare ces normes \`a d'autres plus
traditionnelles, telles que les normes $\|f\|_{\lambda,\mu,\beta}$ de
\eqref{norm1}.

Les normes $\cZ$ poss\`edent des propri\'et\'es remarquables
vis-\`a-vis de la composition et du produit.  Le param\`etre $\tau$
permet de compenser, dans une certaine mesure, la filamentation.
Enfin le caract\`ere hybride est bien adapt\'e \`a la g\'eom\'etrie du
probl\`eme. Si $f$ ne d\'epend que de $x$, la norme \eqref{Z}
co\"{\i}ncide avec la norme d'alg\`bere $\cF^{\lambda\tau+\mu,\gamma}$ d\'efinie par 
\begeq\label{F} \|f\|_{\cF^{\lambda\tau+\mu,\gamma}} =
\sum_{k\in\Z^d} |\hat{f}(k)|\,e^{2\pi(\lambda\tau+\mu)|k|}\,
(1+|k|)^\gamma.
\endeq
(Nous utilisons \'egalement la version \og homog\`ene\fg{}
$\dot{\cF}^{\lambda\tau+\mu,\gamma}$ lorsque l'on n'inclut le mode
$k=0$ dans la somme.)

L'\'equation de Vlasov est r\'esolue par un {\bf sch\'ema de Newton} dont la
premi\`ere \'etape est la solution du lin\'earis\'e autour de $f^0$~:
\[ f^n = f^0 + h^1 + \ldots + h^n, \]
\[
\begin{cases} \pa_t h^1 + v\cdot\nabla_x h^1 
+ F[h^1]\cdot\nabla_v f^0 =0 \\[2mm]
h^1(0,\,\cdot\,) = f_i -f^0
\end{cases} \]
\[
n \ge 1, \quad 
\begin{cases} \pa_t h^{n+1} + v\cdot\nabla_x h^{n+1} + F[f^{n}]\cdot\nabla_v h^{n+1} + F[h^{n+1}]\cdot\nabla_v f^{n} =
- F[h^{n}]\cdot\nabla_v h^{n}\\[2mm]
h^{n+1}(0,\,\cdot\,) = 0.
\end{cases} \]

Dans un premier temps, on \'etablit la r\'egularit\'e analytique en temps petit des $h^n(\tau,\,\cdot\,)$ 
en norme $\cZ^{\lambda,(\mu,\gamma);1}_\tau$; cette \'etape, dans l'esprit d'un th\'eor\`eme de Cauchy--Kowalevskaya,
est r\'ealis\'ee gr\^ace \`a l'identit\'e
\begeq\label{leftddtnorm} \left.\frac{d}{dt}^+ \right|_{t=\tau}
\|f\|_{\cZ^{\lambda(t),\mu(t);p}_\tau} \leq -\frac{K}{1+\tau}\, \|\nabla
f\|_{\cZ^{\lambda(\tau),\mu(\tau);p}_\tau},
\endeq
o\`u $\lambda(t) = \lambda -Kt$, $\mu(t) = \mu - Kt$.
\sm

Dans un deuxi\`eme temps, on \'etablit des estimations uniformes en
temps sur chaque $h^n$, cette fois par une m\'ethode partiellement
eul\'erienne et partiellement lagrangienne, en int\'egrant
l'\'equation le long des caract\'eristiques
$(X^n_{\tau,t},V^n_{\tau,t})$ cr\'ees par la force $F[f^n]$.  (Ici
$\tau$ est le temps initial, $t$ le temps courant, $(x,v)$ les
conditions initiales, $(X^n,V^n)$ l'\'etat courant.) La r\'egularit\'e
de ces caract\'eristiques est exprim\'ee par des contr\^oles en norme
hybride sur les op\'erateurs $\Om^n_{t,\tau}(x,v) =
(X^n_{t,\tau},V^n_{t,\tau})(x+v(t-\tau),v)$, qui comparent la
dynamique perturb\'ee \`a la dynamique libre, informellement
appel\'es {\bf op\'erateurs de scattering} (en temps fini).

On propage alors un certain nombre d'estimations \`a travers le
sch\'ema, les plus importantes \'etant (en simplifiant l\'eg\`erement)
\begeq\label{suphn} 
\sup_{\tau\geq 0}\ \left\| \int_{\R^d} h^n
  \bigl(\tau,\,\cdot\,,v\bigr)\,dv \right\|_{\cF^{\lambda_n\tau +
    \mu_n}}\leq \delta_n,
\endeq
\begeq\label{supttauhn}
\sup_{t\geq\tau\geq 0}\ \Bigl\|h^n \bigl(\tau,\Om^n_{t,\tau}\bigr) \Bigr\|_{\cZ^{\lambda_n(1+b),\mu_n;1}_{\tau-\frac{bt}{1+b}}}
\leq \delta_n,\qquad b= b(t) = \frac{B}{1+t}
\endeq
\begeq\label{Omn}
\Bigl\|\Om^n_{t,\tau}-\Id \Bigr\|_{\cZ^{\lambda_n(1+b),(\mu_n,\gamma);\infty}_{\tau-\frac{bt}{1+b}}}
\leq C\, \left(\sum_{k=1}^n 
\frac{\delta_k\,e^{-2\pi(\lambda_k-\lambda_{n+1}) \tau}}
{2\pi (\lambda_k-\lambda_{n+1})^2}
\right)\,\min\{ t-\tau \, ; \, 1\}.
\endeq 

On note, dans \eqref{suphn}, l'augmentation lin\'eaire de la
r\'egularit\'e de la densit\'e spatiale, qui fait contrepoids
\`a la d\'et\'erioration de r\'egularit\'e dans la variable de
vitesse.  Dans \eqref{supttauhn}, le l\'eger d\'ecalage des indices
par la fonction $b(t)$ sera crucial pour absorber les termes d'erreur
provenant de la composition; la constante $B$ est elle-m\^eme choisie
en fonction des estimations en temps petit r\'ealis\'ees
pr\'ec\'edemment. Enfin, dans \eqref{Omn}, remarquons le contr\^ole
uniforme en $t$, et l'am\'elioration des estim\'ees dans les deux r\'egimes
asymptotiques $t\to\tau$ et $\tau\to\infty$; ceci \'egalement est important
pour la gestion des termes d'erreur. Les constantes $\lambda_n$ et
$\mu_n$ d\'ecroissent \`a chaque \'etape du sch\'ema, convergeant ---
pas trop rapidement --- vers des limites $\lambda_\infty$,
$\mu_\infty$ positives~; dans le m\^eme temps, les constantes $\delta_n$
convergent extr\^emement vite vers $0$, ce qui garantit \og par
r\'etroaction\fg{} l'uniformit\'e des constantes du membre de droite de
\eqref{Omn}.

Les estimations \eqref{Omn} sont obtenues par des applications
r\'ep\'et\'ees de th\'eor\`emes de point fixe en normes
analytiques. Un autre ingr\'edient essentiel pour passer de l'\'etape
$n$ \`a l'\'etape $n+1$ est le m\'ecanisme d'{\bf extorsion de
  r\'egularit\'e}, que nous allons d\'ecrire dans une version
simplifi\'ee.  \'Etant donn\'ees deux fonctions de distribution $f$ et
$\ov{f}$, d\'ependant de $t,x,v$, d\'efinissons
\[ \sigma(t,x) = \int_0^t \int_{\R^d} \bigl( F[f]\cdot\nabla_v
\ov{f}\bigr) \bigl(\tau,x-v(t-\tau),v\bigr)\,dv\,d\tau.\]
L'interpr\'etation de $\sigma$ est comme suit: si les particules
distribu\'ees selon $f$ exercent une force sur des particules
distribu\'ees selon $\ov{f}$, alors $\sigma$ est la variation de
densit\'e $\int f\,dv$ caus\'ee par la r\'eaction de $\ov{f}$ sur
$f$. Nous montrons que si $\ov{f}$ a une r\'egularit\'e glissante
\'elev\'ee, alors la r\'egularit\'e de $\sigma$ en temps grand est
meilleure que ce que l'on attendrait: 
\begeq\label{sigmaleq}
\|\sigma(t,\,\cdot\,)\|_{\dot{\cF}^{\lambda t+\mu}} \leq \int_0^t
K(t,\tau)\,
\bigl\|F\bigl[f(\tau,\,\cdot\,)\bigr]\bigr\|_{\cF^{\lambda\tau+\mu,\gamma}}\,d\tau,
\endeq
o\`u
\[ K(t,\tau) = \left[ \sup_{0\leq s\leq t} \left(\frac{\bigl\|\nabla_v
      \ov{f}(s,\,\cdot\,)\bigr\|_{\cZ^{\ov{\lambda},\ov{\mu};1}_s}}
    {1+s}\right) \right] \,(1+\tau)\, \sup_{k\neq 0,\ \ell\neq 0}\
\frac{e^{-2\pi(\ov{\lambda}-\lambda)|k(t-\tau)+\ell\tau|}\, e^{-2\pi
    (\ov{\mu}-\mu)|\ell|}}{1+|k-\ell|^\gamma}.\] Le noyau $K(t,\tau)$
est d'int\'egrale $O(t)$ quand $t\to\infty$, ce qui laisse craindre
une instabilit\'e violente; mais au fur et \`a mesure que $t$ augmente,
il est \'egalement de plus en plus
concentr\'e sur des temps discrets $\tau = kt/(k-\ell)$; c'est la
manifestation des {\bf \'echos plasmas}, d\'ecouverts et observ\'es
exp\'erimentalement dans les ann\'ees soixante \cite{echo:expe}. Le
r\^ole stabilisant du ph\'enom\`ene d'\'echo, en relation avec
l'amortissement Landau, est mis \`a jour dans notre \'etude.
Le m\'ecanisme d'extorsion de r\'egularit\'e \'eclaire d'un jour nouveau
l'interpr\'etation qualitative populaire de l'amortissement Landau comme
un transfert d'\'energie d'onde vers particule.

On analyse alors la r\'eponse non lin\'eaire due aux \'echos.  Lorsque
$\gamma>11$ on peut montrer que cette r\'eponse est
sous-exponentielle, de sorte que l'on peut la contr\^oler par une
tr\`es l\'eg\`ere perte de r\'egularit\'e glissante, au prix d'une
constante gigantesque, qui sera plus tard absorb\'ee par la
convergence ultrarapide du sch\'ema de Newton. Finalement, une partie
de la regularit\'e de $\ov{f}$ aura \'et\'e convertie en
d\'ecroissance en temps long.

Lorsque $\gamma=1$, des estimations plus fines sont
n\'ecessaires. Pour pouvoir traiter ce cas, nous \'etudions la
r\'eponse non lin\'eaire mode par mode, c'est-\`a-dire en estimant
$\hat \rho(t,k)$ pour tous les $k$, {\it via} un syst\`eme infini
d'in\'egalit\'es. Cela permet de tirer avantage du fait que les
\'echoes qui se produisent \`a des fr\'equences diff\'erentes sont
asymptotiquement plut\^ot bien s\'epar\'es. Par exemple en dimension
$1$, l'\'echo dominant au temps $t$ et \`a la fr\'equence $k$
correspond \`a $\tau = kt/(k+1)$.

En pratique, les trajectoires rectilignes dans \eqref{sigmaleq}
doivent \^etre remplac\'ees par les caract\'eristiques (ceci traduit
le fait que $\ov{f}$ exerce aussi une force sur $f$), ce qui est
source de difficult\'es techniques consid\'erables. Parmi les moyens
mis en oeuvre pour les surmonter, mentionnons un deuxi\`eme
m\'ecanisme d'extorsion de r\'egularit\'e, en temps court et proche
dans l'esprit des lemmes de moyenne; en voici une version
simplifi\'ee: \begeq\label{sigmaleq2}
\|\sigma(t,\,\cdot\,)\|_{\dot{\cF}^{\lambda t+\mu}} \leq \int_0^t
\bigl\|F\bigl[f(\tau,\,\cdot\,)\bigr]\bigr\| _{\cF^{\lambda [\tau -
    b(t-\tau)] +\mu,\gamma}} \, \bigl\|\nabla f(\tau,\,\cdot\,)\bigr\|
_{\cZ^{\lambda(1+b), (\mu,0);1} _{\tau - bt/(1+b)}} \,d\tau.
\endeq
On voit dans \eqref{sigmaleq2} que la r\'egularit\'e de $\sigma$ est
meilleure que celle de $F[f]$, avec un gain qui d\'eg\'en\`ere en
temps grand et lorsque $\tau \to t$. \qed\med





\bibliographystyle{acm}
\bibliography{./note}


\signcm
\signcv

\end{document}